# The upper bound on number of graphs, with fixed number of vertices, that vertices can be colored with n colors.

KAMIL KULESZA, ZBIGNIEW KOTULSKI
Institute of Fundamental Technological Research, Polish Academy of Sciences
ul.Świętokrzyska 21, 00-049, Warsaw Poland, e-mail: {kkulesza, zkotulsk}@ippt.gov.pl

**Abstract:**
In the paper we state and prove theorem describing the upper bound on number of the graphs that have fixed number of vertices $|V|$ and can be colored with the fixed number of $n$ colors. The bound relates both numbers using power of 2, while the exponent is the difference between $|V|$ and $n$. We also state three conjectures on the number of graphs that have fixed number of vertices $|V|$ and chromatic number $n$.

**Keywords:** graph theory, vertex colouring, discrete mathematics, combinatorial problems

## 1. Preliminaries

Graph vertex coloring is an active field of research, with many interesting problems (e.g., [2], [10]). Often these problems belong to NP class of computational complexity, see [4]. As the result, apart from classic approach, they are handled with the help of probabilistic methods (e.g., [6]). Such methods often require answering questions about the bounds, which also motivated our research.

It is widely believed that research into the graph enumeration was initiated in XIX century by Arthur Cayley [1]. The field has proven to be active and fertile, the complete survey was given by F. Harary and E. M. Palmer in [3]. The most powerful tool to count number of non-isomorphic graphs is Pólya's Theorem [7]. To enumerate labeled graphs, the Read's formula can be applied (see [8]), while more general case requires application of Robinson's method, see [9]. Unfortunately, the last one works only for some special situations, while, in general, problem remains unsolved, see [3].

The Pólya's Theorem, which is dealing with non-isomorphic graphs, at least to some extend, was motivated by chemical enumeration (see [7]). In our research graph coloring was used to construct check-digit scheme, see [5]. Thus, different questions had to be asked. The answers required taking into account some isomorphic graphs.

The paper deals with the vertex coloring of graphs that are finite and simple. The graph coloring considered is proper vertex $n$-coloring ($n \geq \chi(G)$). It does not need to be the minimal coloring ($n = \chi(G)$), unless stated otherwise. The main goal of the paper is to state and prove the Upper Bound Theorem.

*Notation*

Let $\Gamma(|V|)$ be the number of all possible graphs having $|V|$ vertices. It is well known that fact $\Gamma(|V|) = 2^{\frac{|V|(|V|-1)}{2}}$. In other words, $\Gamma(|V|)$ represents number of binary sequences such that every bit corresponds to the possible edge in the $|V|$ vertices graph (1 denotes presence of the edge, 0 otherwise). Let's denote:

$x_j$ as the number of vertices colored by $j$th color $(0 < j \leq n)$. Vertices are partitioned into $n$ sets, such that all vertices with the same color are in one set;

$V_j$ as the set of vertices colored by $j$th color $(0 < j \leq n)$, $|V_j| = x_j$;

$P$ as particular, no-degenerate partition of $|V|$ into $n$ color sets. $P = \{x_1, x_2, ..., x_j, ..., x_n\}$ and $\sum_{j=1}^{n} x_j = |V|$. All vertices are homogenous, the only difference among them is the color ;

$\Gamma(|V|, n, P)$ as the number of the graphs having $|V|$ vertices and particular $P$ ($|V|$ partition into $n$ pieces). Applying reasoning analogous to one concerning $\Gamma(|V|)$ one can observe that

$$\Gamma(|V|, n, P) = 2^{(x_2 + x_3 + ... + x_n)x_1 + (x_3 + x_4 + ... + x_n)x_2 + ... + x_n x_{n-1}} \quad (1)$$

$\Gamma(|V|, n) = \max_P \{\Gamma(|V|, n, P)\}$ where $\max_P$ denotes maximum value of $\Gamma(|V|, n, P)$ over all possible partitions $P$ for given $|V|$ and $n$. In plain English $\Gamma(|V|, n)$ represents such partition $P$, that yields maximum number of graphs, while both $|V|$, $n$ are fixed.

## 2. Main result

Having notation in place, it is time to state and prove the theorem on the upper bound of $\Gamma(|V|, n)$.

**Theorem 1 (Upper bound)**
Let $|V| = n + y$ where $y \in N$. Then

$$\Gamma(|V|) \geq 2^y \Gamma(|V|, n) \quad (2)$$

**Proof of Theorem 1.**
*Observation*
Recall that :

a. $\Gamma(|V|) = 2^{\frac{|V|(|V|-1)}{2}}$

b. $\Gamma(|V|, n) = \max_P (\Gamma(|V|, n, P)) = 2^{\max_P \{(x_2 + x_3 + ... + x_n)x_1 + (x_3 + x_4 + ... + x_n)x_2 + ... + x_n x_{n-1}\}}$

Rewriting Theorem 1 yields:

$2^{\frac{|V|(|V|-1)}{2}} = \Gamma(|V|) \geq 2^y \Gamma(|V|, n) = 2^{y + \max_P \{(x_2 + x_3 + ... + x_n)x_1 + (x_3 + x_4 + ... + x_n)x_2 + ... + x_n x_{n-1}\}}$

thus,

$$2^{\frac{|V|(|V|-1)}{2}} \geq 2^{y + \max_P \{(x_2 + x_3 + ... + x_n)x_1 + (x_3 + x_4 + ... + x_n)x_2 + ... + x_n x_{n-1}\}} \quad (3)$$

Equation 3 shows that in order to prove the theorem one needs to compare the exponents. Hence, further work will be carried out in base 2 logarithms only.
Observe:
$\log_2 \Gamma(|V|) = (|V| - 1) + (|V| - 2) + ... + (|V| - i) + ... + (|V| - |V| + 1)$, where $(|V| - i)$ is $i$th term of the sequence;
$\log_2 \Gamma(|V|, n, P) = (x_2 + x_3 + ... + x_n)x_1 + (x_3 + x_4 + ... + x_n)x_2 + ... + x_n x_{n-1}$;
$\log_2 \Gamma(|V|, n) = \max_P \{\log_2 \Gamma(|V|, n, P)\}$, where $\max_P$ denotes maximum over all possible partitions $P$ for given $|V|$ and $n$.

Proving the theorem requires finding partitions $P$ that yield maximum value of $\log_2 \Gamma(|V|, n, P)$. Two cases will be considered: trivial case with single partition and the case of multiple partitions.

***Case 1:*** $|V| = n \Rightarrow y = 0$

There is only one possible partition $P = \{1,1,...,1\}$ and $\forall j\ x_j = 1$.

It means that in this case $\log_2 \Gamma(|V|, n) = \log_2 \Gamma(|V|, n, P)$

Then $\log_2 \Gamma(|V|, n, P) = (n-1)*1 + (n-2)*1 + ... + 1*1 = \dfrac{(n-1)n}{2}$, because there are *n-1* terms in the summation. But $|V| = n$, hence $\log_2 \Gamma(|V|, n, P) = \dfrac{(|V|-1)|V|}{2}$. Finally:

$$\log_2 \Gamma(|V|, n, P) = \log_2 \Gamma(|V|, n) = \log_2 \Gamma(|V|) \quad (4)$$

yields

$$\Gamma(|V|) = 2^0 \Gamma(|V|, n) = 2^y \Gamma(|V|, n) \quad (4a)$$

Equation 4a completes proof for the Case 1.

***Case 2:*** $|V| > n \Rightarrow y > 0$

*Observation:*

In this case $\log_2 \Gamma(|V|, n, P) = (x_2 + x_3 + ... + x_n)x_1 + (x_3 + x_4 + ... + x_n)x_2 + .. + x_n x_{n-1}$ has $n-1$ terms. Due to the fact that $\sum_{j=1}^{n} x_j = |V|$, $\forall P\ \exists j_0$ such that $x_{j_0} > 1$. Term from the expansion of $\log_2 \Gamma(|V|, n, P)$ with $x_{j_0}$ as the multiplier will have the from: $(x_{j_0+1} + x_{j_0+2} + ... + x_n)x_{j_0}$, which can be expressed as sum of $x_{j_0}$ identical terms $(x_{j_0+1} + x_{j_0+2} + ... + x_n)$. When such expansion is used for all $j_0$ ($x_{j_0} > 1$), then $\log_2 \Gamma(|V|, n, P)$ has $(|V|-1)$ terms as $\log_2 \Gamma(|V|)$.

This observation allows to compare both sequences term-wise. Convenient way to do it is to enumerate all terms of $\log_2 \Gamma(|V|, n, P)$ in the same way as $\log_2 \Gamma(|V|)$ terms. To establish the correspondence between *i*- indexing of $\log_2 \Gamma(|V|)$ terms and original (before expansion) *j*-indexing of $\log_2 \Gamma(|V|, n, P)$, define:

$i$ is $\log_2 \Gamma(|V|)$ terms' enumeration, it is also $\log_2 \Gamma(|V|, n, P)$ terms' enumeration when $|V| = n$ (Case 1), $i = x_1 + x_2 + ... + x_j = \sum_{q=1}^{j} x_q = |V|$

$j$ is $\log_2 \Gamma(|V|, n, P)$ terms' enumeration when $|V| > n$

$S_1^i$ as the sum of *ith* term elements in $\log_2 \Gamma(|V|)$. $S_1^i = (x_{i+1} + x_{i+2} + ... + x_n)$.

$S_2^i$ as the sum of *ith* term elements in $\log_2 \Gamma(|V|, n, P)$. Formula for $S_2^i$, relation between $S_2^i$ and $S_1^i$ have to be found. Further argument will be illustrated by the example.

***Example1*** (inside the proof).

Take $|V| > n$ such that $y \geq 3$ and some $P = \{1,3,1,...\}$. Then:

$\log_2 \Gamma(|V|) = (|V|-1) + (|V|-2) + .. + (|V|-i) + .. + (|V|-|V|+1) = S_1^1 + S_1^2 + S_1^3 + S_1^4 + .. + S_1^i + .. + S_1^{|V|-1}$ (5)

$\log_2 \Gamma(|V|, n, P) = (x_2 + x_3 + ... + x_n)*1 + (x_3 + x_4 + ... + x_n)*3 + (x_4 + x_5 + ... + x_n)*1 + ... =$
$= (|V|-1)*1 + (|V|-4)*3 + (|V|-5)*1 + ...$

Expanding term $(x_3 + x_4 + ... + x_n)*3$ into addition yields :

$$\log_2 \Gamma(|V|, n, P) = (|V|-1) + (|V|-4) + (|V|-4) + (|V|-4) + (|V|-5) + \ldots \qquad (6)$$

Writing equation 6 as the sum of $S_2^i$ terms:

$$\log \Gamma_2(|V|, n, P) = S_2^1 + S_2^2 + S_2^3 + S_2^4 + S_2^5 + \ldots \qquad (7)$$

Using equations 6 and 5 observe that $S_2^2 = S_2^3 = S_2^4 = S_1^4$, so equation 7 can be written as:

$$\log_2 \Gamma(|V|, n, P) = S_2^1 + S_2^2 + S_2^3 + S_2^4 + S_2^5 + \ldots = S_1^1 + S_1^4 + S_1^4 + S_1^4 + S_1^5 + \ldots \qquad (8)$$

*End of example.*

In general case, when $\log_2 \Gamma(|V|, n, P)$ is expanded into $(|V|-1)$ terms the following holds: $S_1^i \geq S_2^i$ for all corresponding sums of terms. It yields that $\log_2 \Gamma(|V|) \geq \log_2 \Gamma(|V|, n, P)$.

When $|V| > n$, $\exists i$ such that $S_1^i > S_2^i$ (it results from at least one $x_{j_0} > 1$). It will be shown, that this property yields $\log_2 \Gamma(|V|) > \log_2 \Gamma(|V|, n, P)$.

Define:

$\Psi(S_2^{j_0}) = (x_{j_0+1} + x_{j_0+2} + \ldots + x_n)x_{j_0} = S_1^{j_0} * x_{j_0}$ is the sum of $x_{j_0}$ identical terms of $\log_2 \Gamma(|V|, n, P)$ resulting from $x_{j_0} > 1$. In the example above: $\Psi(S_2^2) = S_1^4 + S_1^4 + S_1^4$;

$\Psi(S_1^{i_0}) = S_1^{(i_0 - x_{j_0})+1} + S_1^{(i_0 - x_{j_0})+2} + \ldots + S_1^{i_0-1} + S_1^{i_0}$ is the sum of $x_{j_0}$ terms of $\log_2 \Gamma(|V|)$ corresponding to $x_{j_0}$ identical terms of $\log_2 \Gamma(|V|, n, P)$ ($i_0$ is derived from $j_0$). In the example above $\Psi(S_2^4) = S_1^2 + S_1^3 + S_1^4$;

Observe that:

$$\Psi(S_1^{i_0}) = S_1^{i_0 - x_{j_0}+1} + S_1^{i_0 - x_{j_0}+2} + \ldots S_1^{i_0-1} + S_1^{i_0} = \frac{1}{2}\left(S_1^{i_0+1-x_{j_0}} + S_1^{i_0}\right)x_{j_0} \qquad (9)$$

Since $S_1^{i_0+1-x_{j_0}} = S_1^{i_0} + x_{j_0} - 1$, hence

$$\Psi(S_1^{i_0}) = \frac{1}{2}\left(S_1^{i_0+1-x_{j_0}} + S_1^{i_0}\right)x_{j_0} = \frac{1}{2}\left(2 S_1^{i_0} + x_{j_0} - 1\right)x_{j_0} = S_1^{i_0} x_{j_0} + \frac{1}{2}(x_{j_0} - 1)x_{j_0} \qquad (10)$$

On the other hand $\Psi(S_2^{j_0}) = (x_{j_0+1} + x_{j_0+2} + \ldots + x_n)x_{j_0} = S_2^{i_0} x_{j_0} = S_1^{i_0} x_{j_0}$, because $S_2^{i_0} = S_1^{i_0}$. This observation allows to state the relation between $\Psi(S_1^{i_0})$ and $\Psi(S_2^{j_0})$:

$$\Psi(S_1^{i_0}) = \Psi(S_2^{j_0}) + \frac{1}{2}(x_{j_0} - 1)x_{j_0} \qquad (11)$$

for $x_{j_0} > 1$, yields

$$\Psi(S_1^{i_0}) > \Psi(S_2^{j_0}) \qquad (12)$$

This shows that

$$\log_2 \Gamma(|V|) > \log_2 \Gamma(|V|, n, P) \qquad (13)$$

provided it exits $j_0$ such that $x_{j_0} > 1$. Hence equation 13 holds for all $n$ such $n < |V|$ (since $n < |V| \Rightarrow \exists j_0$ such that $x_{j_0} > 1$).

Recall that $\log_2 \Gamma(|V|, n) = \max_P \{\log_2 \Gamma(|V|, n, P)\}$, so for $n < |V|$ it holds that

$$\log_2 \Gamma(|V|) > \log_2 \Gamma(|V|, n) \qquad (14)$$

Having equation 14, it is time to show where does power of 2 in the equation 2 comes from. We start our reasoning from the point where

$$\lambda = \log_2 \Gamma(|V|) - \log_2 \Gamma(|V|, n) \qquad (15)$$

attains minimal value. Using equation 11 one can observe that

$$\lambda = \Psi(S_1^{i_0}) - \Psi(S_2^{j_0}) = \frac{1}{2}(x_{j_0} - 1)x_{j_0} \tag{16}$$

To find minimal value of $\lambda$, the expression $\frac{1}{2}(x_{j_0} - 1)x_{j_0}$ should be minimized for each $x_{j_0} > 1$ in the partition $P$. For $x_{j_0} > 1$, we get $\min\left\{\frac{1}{2}(x_{j_0} - 1)x_{j_0}\right\} = 1$ when $x_{j_0} = 2$. In such case $\Psi(S_1^{i_0}) = \Psi(S_2^{j_0}) + 1$.

If given partition $P$, for each $x_{j_0} > 1$ holds that $x_{j_0} = 2$, then the number of all different indexes $j_0$ is $|V| - n = y$. This yields

$$\log_2 \Gamma(|V|) = \log_2 \Gamma(|V|, n) + y \tag{17}$$

results in

$$\Gamma(|V|) = 2^y \Gamma(|V|, n) \tag{17a}$$

Now we discuss the situation when $|V|$ cannot be partitioned into $n$ pieces such $\forall j \; x_j \leq 2$. For the given $\log_2 \Gamma(|V|, n) \; \exists j_0$ such that $x_{j_0} > 2$ expression $\frac{1}{2}(x_{j_0} - 1)x_{j_0} \geq 3$ and it increases faster then $y$. In this case

$$\log_2 \Gamma(|V|) > \log_2 \Gamma(|V|, n) + y \tag{18}$$

results in

$$\Gamma(|V|) > 2^y \Gamma(|V|, n) \tag{18a}$$

Equation 18a completes proof for the Case 2.

Combining equations 4a, 17a, 18a yields:

$$\Gamma(|V|) \geq 2^y \Gamma(|V|, n) \tag{19}$$

**Q.E.D.**

**Corollary 1**
Let $|V| = n + y$ where $y \in N$. For $2n \geq |V|$ it holds that

$$\Gamma(|V|) = 2^y \Gamma(|V|, n) \tag{20}$$

**Proof**
Condition $2n \geq |V|$ allows to building the partition $P$ such that $\forall j \; x_j \leq 2$. Hence the result follows from argument leading to equations 17 and 17a.
**Q.E.D**

## 3. On the number of graphs such $\chi(G) = n$

The Upper Bound theorem discusses number of graphs that can be colored with $n$ colors. The result includes whole class of graphs $G'$ such that $\chi(G') < n$. This observation can be used to calculate the number of the graphs that have fixed chromatic number $\chi(G)$ (say $\chi(G) = n$).

*Let us denote:*
$\Gamma(|V|, \chi(G) = n, P)$ as the number of the graphs with $\chi(G) = n$, that have $|V|$ vertices and particular $P$ ($|V|$ partition into $n$ pieces);

$\Gamma(|V|, \chi(G) = n) = \max_P \{\Gamma(|V|, \chi(G) = n, P)\}$ where $\max_P$ denotes maximum value of $\Gamma(|V|, \chi(G) = n, P)$ over of all possible partitions $P$ for given $|V|$ and $\chi(G) = n$. In plain English $\Gamma(|V|, \chi(G) = n)$ represents partition $P$, resulting in maximum number of graphs, that cannot be colored with less then $n$ colors.

Having the definitions in place, three conjectures are stated:

**Conjecture 1**

$$\Gamma(|V|, \chi(G) = n) = \Gamma(|V|, n) - \Gamma(|V|, n-1) \tag{21}$$

**Conjecture 2**

$$\Gamma(|V|, n) \geq 2 * \Gamma(|V|, n-1) \tag{22}$$

Combining two conjectures (eq. 21,22) allows us to state final hypothesis in this section.

**Conjecture 3**

$$\Gamma(|V|, \chi(G) = n) \geq \Gamma(|V|, n-1) \tag{23}$$

## 4. Further research and concluding remarks

We stated and proven the Upper Bound Theorem together with resulting Corollary 1. Following these results three conjectures were stated. Our research was motivated by practical considerations resulting from construction of the graph coloring check-digit scheme, see [5]. Obviously, theoretical questions are not limited to the conjectures stated above, two more follows as example.

1. Bounds for more specific conditions.

Using Theorem 1 and it's proof one can derive results stating the relations between $\Gamma(|V|)$ and $\Gamma(|V|, n)$ for more specific conditions (e.g., $\exists j_0$ such that $x_{j_0} = k$, $k \geq 2$).

2. Non-isomorphic graphs.

In our calculation we count number of the graphs for every partition. This leads to the situation that some isomorphic graphs are counted more than once. Obtaining result for non-isomorphic graphs should allow to lower the upper bound in some cases.